\documentclass{amsart}

\usepackage{amssymb}
\usepackage[all]{xy}
\usepackage{hyperref}

\newtheorem{thm}{Theorem}

\newtheorem{lem}[thm]{Lemma}
\newtheorem{cor}[thm]{Corollary}

\newtheorem{prop}[thm]{Proposition}

\theoremstyle{definition}
\newtheorem{defn}[thm]{Definition}

\newtheorem{say}[thm]{}
\newtheorem{exmp}[thm]{Example}

\newtheorem*{ack}{Acknowledgments}      

\newtheorem{defn-thm}[thm]{Definition--Theorem}  
\newtheorem{defn-lem}[thm]{Definition--Lemma}  

\newtheorem{surprise}[thm]{Surprise}

\theoremstyle{remark}

\renewcommand{\c}[0]{{\mathbb C}}  

\renewcommand{\o}[0]{{\mathcal O}} 
\newcommand{\z}[0]{{\mathbb Z}}
\newcommand{\n}[0]{{\mathbb N}}
\renewcommand{\r}[0]{{\mathbb R}} 

\renewcommand{\a}[0]{{\mathbb A}}

\newcommand{\p}[0]{{\mathbb P}}
\newcommand{\f}[0]{{\mathbb F}}
\newcommand{\q}[0]{{\mathbb Q}}

\newcommand{\qtq}[1]{\quad\mbox{#1}\quad}
\newcommand{\spec}[0]{\operatorname{Spec}}

\newcommand{\rank}[0]{\operatorname{rank}}

\newcommand{\proj}[0]{\operatorname{Proj}}

\newcommand{\sing}[0]{\operatorname{Sing}}

\newcommand{\chow}[0]{\operatorname{Chow}}

\newcommand{\hilb}[0]{\operatorname{Hilb}}
\newcommand{\grass}[0]{\operatorname{Grass}}

\newcommand{\tsum}[0]{\textstyle{\sum}}

\def\into{\DOTSB\lhook\joinrel\to}

\def\loccoh#1.#2.#3.#4.{H^{#1}_{#2}(#3,#4)}

\DeclareMathAlphabet{\mathchanc}{OT1}{pzc}%
                                {m}{it}

\newcommand{\GL}{\mathrm{GL}}
\newcommand{\PGL}{\mathrm{PGL}}
\newcommand{\PSL}{\mathrm{PSL}}

\newcommand{\SL}{\mathrm{SL}}

\newcommand{\Sp}{\mathrm{Sp}}

\newcommand{\sym}[0]{\operatorname{Sym}}

\usepackage[all]{xy}\xyoption{dvips}

\begin{document}
\bibliographystyle{amsalpha}

 \title[Mumford and the of  moduli varieties]{Mumford's influence  on the \\ moduli theory of algebraic varieties}
 \author{J\'anos Koll\'ar}

\begin{abstract} We give a short appreciation of Mumford's work on the moduli of varieties by putting it into historical context. By reviewing earlier works we highlight  the  innovations introduced by Mumford. Then we discuss  recent developments whose origins can be traced back to Mumford's  ideas.
\end{abstract}

 \maketitle

  The theory of moduli of algebraic varieties aims to solve the following problems.
\medskip

$\bullet$ (Main question, old style)
       Can we parametrize the set of all varieties of a given class  in a natural way using the points of another algebraic variety?

 \medskip

$\bullet$ (Main questions, new style)
        What is a ``good family'' of algebraic varieties?
  Can we describe all ``good families'' in an  ``optimal'' manner?
\medskip

The old style question was first studied by Riemann \cite{MR1579035} for
algebraic curves---equivalently, compact Riemann surfaces---and most related work for the next century focused on this problem. 

The shift to the new style question happened with Grothendieck's lectures at the Cartan seminar in 1960/61. At the same time  Mumford developed a general framework to approach moduli problems in algebraic geometry---called Geometric Invariant Theory, usually abbreviated as  GIT---and used it to  
complete the construction of $M_g$.

We aim to highlight Mumford's pivotal role in moduli theory, 
how his ideas were different from  earlier work and how they persist in the  subsequent  developments.

\begin{ack} These notes are an expanded version of my lecture at the conference
{\it From Algebraic Geometry to Vision and AI: A Symposium Celebrating the Mathematical Work of David Mumford,}  organized by the 
Center of Mathematical Sciences and Applications under the direction of
Shing-Tung Yau. 

I thank  Rahul~Pandharipande  and Chenyang~Xu for helpful  comments and references. 
Partial  financial support    was provided  by  the NSF under grant number
 DMS-1362960.
\end{ack}

\section{Moduli of curves, 1857--1960}\label{sec.1}

An algebraic curve can be studied using its projective embeddings and their equations, while a compact Riemann surface of genus $\geq 2$ can be viewed as a quotient of the complex disc by a  discrete subgroup of $\PSL_2(\r)$. 
The former leads to an algebraic theory of moduli, the latter to an analytic or group theoretic version. While these 2 approaches frequently intersect, it is convenient to think of them as parallel developments.

\subsection*{Moduli of curves, algebraic theory}{\ }

It has been known at least since Descartes  that the set of all plane curves of degree $d$ 
is naturally parametrized by the vector space ${\mathbf V}_d$ of all polynomials
$f(x,y)$ of degree $\leq d$; later this was replaced by the 
projective space ${\mathbf P}_d$ of all homogeneous polynomials
$f(x,y,z)$ of degree $= d$. (Thus ${\mathbf P}_d$ has dimension
$\binom{d+2}{2}-1$.) It was also clear that curves obtained by a linear change of coordinates are the ``same'' and thus the orbit spaces   ${\mathbf V}_d/\GL_3$ or  ${\mathbf P}_d/\PGL_3$ are the essential objects.

Classical invariant theory was  initiated by Cayley in \cite{cayley45}, building on earlier work of Boole \cite{boole41}. Its aim was to start with a 
vector space ${\mathbf V}$ acted upon by  a group $G$, and then describe all
  $G$-invariant polynomial or rational functions on ${\mathbf V}$.
Initially the main interest was in concrete cases.


If the  orbit space ${\mathbf P}_d/\PGL_3$ were an algebraic variety, then  we could identify rational functions on ${\mathbf P}_d/\PGL_3$  
with $\PGL_3$-invariant  
 rational functions on  ${\mathbf P}_d$. Although the  orbit space is not an algebraic variety,
 the field of $\PGL_3$-invariant rational functions is the function field of a variety---called Mumford's GIT quotient and  denoted by  ${\mathbf P}_d/\!/\PGL_3$---that serves as the ``best algebraic approximation'' of the  orbit space   ${\mathbf P}_d/\PGL_3$. 
However,  the classical authors do not seem to have focused on this aspect of invariant theory. Their main interest was in getting explicit generators for the ring of invariants in the affine case ${\mathbf V}_d/\GL_3$ and its generalizations. Eventually   Gordan proved \cite{gordan68}  that
rings of invariants of binary forms are finitely generated and 
Hilbert   proved  finite generation for any reductive group action  \cite{MR1510634}. However, the original hope of Cayley and Gordan of getting explicit generators is still out of reach.


Another moduli-theoretic direction was also initiated by
Cayley, who, in  two articles with the same title,
 constructed  the moduli space  of space curves
\cite{cayley1860, cayley1862}, which  is now usually called the {\it Chow variety} of curves in $\p^3$.

To any curve $C\subset \p^3$, Cayley's method associates the set of all lines  $L \subset \p^3$  that meet $C$. This set is a hypersurface
in the Grassmannian of all lines $\grass(1,3)$,  giving an injection
$$
\left\{
\mbox{degree $d$ curves in $\p^3$}
\right\}
\into
 \bigl| \o_{\grass(1,3)}(d)\bigr|
$$
whose image is an algebraic variety. 
The general correspondence between
$m$-cycles in  $\p^N$ and hypersurfaces in  $\grass(N-m-1,N)$ 
was fully worked out only in   \cite{MR1513117}.  
A complete algebraic treatment  is given in \cite[Secs.X.6--8]{hodge-ped};
see \cite[Sec.I.3]{rc-book} for a more geometric one.


Severi 
proved that  $M_g$ (the moduli space of smooth projective curves of genus $g$) is unirational for $g\leq 10$ \cite{severi-15, MR0245574}, but it is not clear to me what
 Severi thought about the existence of $M_g$ as an algebraic variety. 
Note that one can state and prove that $M_g$ is  unirational without
saying anything about its existence. Severi proves that for certain values
$(g,d,r)$, the set of degree $d$ plane curves with $r$ nodes is a rational variety and a general curve of genus $g$ is birational to a degree $d$ plane curve with $r$ nodes. Even with  modern techniques this approach proves only that some open subset of $M_g$ is an algebraic variety for $g\leq 10$.

It seems to have been Weil who first seriously contemplated that $M_g$ and $A_g$
(the moduli space of principally polarized Abelian varieties) should be  varieties defined over $\q$. This is a quite non-trivial observation since the deepest classical approaches to $M_g$ via $T_g$ (the Teichm\"uller space) or $A_g$
rely on transcendental functions.
If $C$ is a curve of genus $g$ and $[C]\in M_g$ denotes the corresponding point, then its residue field  $k_C:=k\bigl([C]\bigr)$ is an invariant of $C$,
 called the {\it field of moduli} of $C$. A key question is whether
$C$ can be defined over its field of moduli; now we see this as the
distinction between fine and coarse moduli spaces.
Weil and  Matsusaka discussed the fields of moduli in various examples
\cite{MR0023093, MR0094360}.

Satake  \cite{MR0170356} and  Baily-Borel \cite{MR0216035} worked out compactifications of  $A_g$
and of other quotients of  symmetric spaces, but
in their works the boundary exists mostly as a set-theoretic object.
It is very misleading for the general moduli theory that the points of the boundary of $A_g$ naturally correspond to lower dimensional Abelian varieties, precluding the existence of universal families.
This may have been a true psychological barrier historically.

Igusa gave an explicit and  complete description of $M_2$, and what is now called its GIT compactification,  over $\spec \z$, building  on classical invariant theory  \cite{MR0114819}.

The first explicit claim that $M_g$ should exist as an algebraic variety may be in \cite{weil-ea}:
 ``As for $M_g$ there is virtually  no doubt that it can be provided with the structure of  an algebraic variety.''

\subsection*{Moduli of curves, analytic theory}{\ }
      
   A very thorough discussion of this topic is given in the papers by
A'Campo, Ji and Papadopoulos \cite{MR3021551, ajp}. Anyone interested in the details should start by reading these  papers and   the extensive bibliography in \cite[Vol.III]{MR0257326}. Here I just mention the key works quickly to give an overview.

The story starts with 
Riemann  \cite{MR1579035}
who proves that  Riemann surfaces of genus $g$ depend on $3g-3$ parameters. In modern language he explains that the local deformation space of a Riemann surface of genus $g$ has complex dimension $3g-3$. 
Riemann introduces the word ``Mannigfaltigkeiten''  (usually translated as `manifoldness'') to refer to $M_g$, but neither Riemann nor later Klein were clear on what this precisely means. 
(A formal definition of complex manifolds appears in the literature only after 1940 in the works of Teichm\"uller, Chern and Weil.)


 Following Riemann's approach,  Hurwitz  \cite{MR1510692}
defines what we now call the Hurwitz scheme and uses it to prove that 
 $M_g$ is irreducible. 

The approach through Fuchsian groups and automorphic functions is developed in the monumental 1300 page work of   
Fricke and Klein \cite{fri-kle-1, fri-kle-2}.
They prove that   $T_g$ exists and that it is a contractible topological space.
These books contain many examples where a  complex structure is also given, but there does not seem to be a general claim that $T_g$ always has a  ``natural'' complex structure.

The first  modern treatments is due to Siegel 
who  constructed  $A_g$ as an analytic space  \cite{MR1503238}.
The complete theory of $T_g$ is worked out by 
Teichm\"uller  \cite{teichm}, though it took some time before this was fully understood by experts through the works of Ahlfors, Bers and Weil.

It was probably the Bourbaki talks of Weil \cite{MR0124485b, MR0124485} that---though focusing on analysis and differential geometry---brought these developments to the full attention of algebraic geometers.


The shift from the old  to the new style  questions, still in the analytic setting,   was accomplished  by
Grothendieck's lectures \cite{MR0106907b}.
Grothendieck defines the functor of smooth families of compact Riemann surfaces 
(over reduced analytic spaces) and proves that the Teichm\"uller
space   $T_g$ represents this functor. 
Grothendieck  mentions his plan to work out   most of the theory  algebraically over arbitrary  schemes as well. He states that he has other ideas that should construct $M_g^{[r]}$---the moduli of curves
with level structure $r\geq 3$---as a scheme over $\spec \z$. 
(Note that if $r\geq 3$ then $M_g^{[r]}$ is a fine moduli space.)
Grothendieck also says that quasi-projectivity may not follow from his methods. (This is a potentially serious issue. Although $M_g$ is a quotient of $M_g^{[r]}$ by $\Sp_{2g}(\z/r)$, the quotient of a non-quasi-projective variety by a finite group need not be a scheme. Thus one cannot yet  conclude that $M_g$ is a scheme, much less that it is  quasi-projective.)
Grothendieck also refers to Mumford's work that was expected to construct $M_g$  as a quasi-projective scheme over $\spec \z$.


 \section{Mumford's work on the moduli of curves}
     \label{sec.2} 


{\it When [Zariski] spoke the words algebraic variety, there was a certain resonance in his voice that said distinctly that he was looking into a secret garden. I immediately wanted to be able to do this too ... Especially, I became obsessed with a kind of passion flower in this garden, the moduli spaces of Riemann.} (Mumford in \cite[p.225]{mumf-auto}.)

\medskip

This description, reminiscent  of Newton 
``diverting [himself] in now and then finding a smoother pebble or a prettier shell than ordinary,'' may explain the origin of Mumford's interest in $M_g$, but it does not hint at his decades-long systematic work to understand not just $M_g$, but moduli theory in general.

 
Next we list the main developments introduced by Mumford, explain how these were different from earlier methods, and then  the rest of the paper  outlines how they led to a moduli theory of higher dimensional varieties.
It is rather presumptuous to claim to know what 
      Mumford's main ideas were without doing  much more serious historical research. So the following  topics are those that either strike me as having been especially novel,   or that  had  a major  influence on my own studies of moduli.

\begin{say}[Mumford's main ideas] \label{main.ideas}
I will focus on 5 topics.
\begin{enumerate}
\item   $M_g$ coarsely represents the functor of families of smooth, projective  curves of genus $g$.
\item GIT gives the right moduli spaces.
\item   $M_g$  should be constructed as a quotient of the space of 
 $m$-canonically embedded curves   
 $\{C\into \p^{N-1}\}$   by the group $\PGL_N$ for $N=(2m-1)(g-1)$.
\item   $M_g$  should be compactified using reducible nodal curves  (with Deligne).
\item   $M_g$ and $\bar M_g$ are  non-linear analogs of 
 Grassmannians; hence their   cohomology rings  govern many 
 enumerative questions.
\end{enumerate}

It is sometimes hard to appreciate how transformative some  definitions and theorems  were 60 years ago. 
Several of Mumford's ideas  seem to have been quite alien to
his predecessors, but once Mumford introduced them,  they  quickly became viewed as the
``obvious approach.''

As we discussed, thinking of moduli spaces as representing functors did not appear until \cite{MR0106907b}. While the functorial point of view quickly gained acceptance with the expansion of Grothendieck's school, $M_g$ may have seemed an exception since it does not represent the functor of families of smooth, projective  curves of genus $g$. (In some sense we are still trying to come to terms with this using stacks.) Mumford's pragmatic  balance of the ``pure thought'' approach of Grothendieck and the classical authors probably appeared unnatural to both sides at first.

While the rings of invariants of many classical group actions were much studied, there was no  systematic effort to view them as giving compact moduli spaces. A key observation of Mumford is that one can identify the  ``good'' points  (now known as stable and semi-stable points) and obtain a moduli space for the corresponding objects.  
This approach has been spectacularly successful for sheaves; see  \cite{MR0175899} for the early results and 
 \cite[Chap.4]{huy-leh} for a recent survey.
We discuss in  Section~\ref{sec.4} that
we have not yet seen the full story for higher dimensional varieties.


\medskip
{\it Note on terminology.} In the literature, the adjective ``stable'' is used in at least 3 distinct senses. For clarity I will use {\it GIT-stable}
for the geometric invariant theory notion, {\it DM-stable} or
{\it Deligne-Mumford-stable} for stable curves as in Theorem~\ref{dm.comp.thm} and
{\it KSB-stable} or {\it Koll\'ar--Shepherd-Barron-stable} for the
higher dimensional version as in Definition~\ref{KSB.stab.defn}. 
\medskip

I am rather mystified why   (\ref{main.ideas}.3) was not proposed much earlier.
There seems to have been a strong preference for viewing curves either as covers of $\p^1$ or as plane curves (usually with nodes). However people did study curves in $\p^n$. Another problem may have been the---probably correct---feeling that this approach could not be used to understand the geometry of $M_g$ anyhow. In particular, invariant theory was usually considered only for polynomial rings before Mumford.

Classical authors  developed a very good understanding of singular, but usually  irreducible, plane curves. Thus in studying degenerations of curves, there was probably a strong preference for  irreducible limit curves. 

The compactification of $M_1$ is also misleading here. We get a single point at infinity but, depending on how we write our elliptic curves, we get a nodal rational curve or a cycle of rational curves in natural families.

We illustrate the difference between previous approaches and the Deligne-Mumford compactification of genus 2 curves in Section~\ref{sec.3}.

For the higher dimensional theory it was clear that one should understand the correct  generalization of nodes. The problem is that nodes have many good properties and one needs to choose the right one to generalize. We come back to this in Section~\ref{sec.6}.

Even with hindsight I have no idea how Mumford knew that $M_g$ will turn out to be a  ``passion flower.'' The 2018ICM lectures of Pandharipande  give a  very good recent survey about the beauty of $M_g$; I recommend both the video  \cite{rahul-icm-v}   and the written form \cite{rahul-icm-w}.  
\end{say}

The above list covers only a  part of Mumford's work on moduli. Several other  results that turned out to be important for the general theory
are mentioned later, but other topics, like the compactification of $A_g$ are missing \cite{amrt}.
One needs to read the
 complete collection of Mumford's papers and books to get a full appreciation.
These are very conveniently  available at
\url{http://www.dam.brown.edu/people/mumford}.

\section{Example of genus 2 curves}
     \label{sec.3} 
   
We illustrate the difference between the old-style set theoretic approach and the new-style functorial approach with genus 2 curves, where everything can be written down quite explicitly. For simplicity assume that we work over $\c$.

Let  $C$ be a smooth, projective curve of genus 2; equivalently, 
a smooth, compact Riemann surface of genus 2.
Then  $h^0(C, \omega_C)=2$, thus we get
 a unique degree 2 morphism  $\tau:C\to {\mathbb P}^1$
which  ramifies at 6 points. So we can  write  $C$ as   
$$
\bigl(z^2=f_6(x{:}y)\bigr)\subset \p^2(1,1,3),
$$  
where $f_6$ is a homogeneous polynomial of degree 6 and
the notation $\p^2(1,1,3)$ means that we declare the $z$ variable to have degree 3. In particular, the usual homegeneity relation changes to $(x{:}y{:}z)= \bigl(\lambda x{:}\lambda y{:}\lambda^3z\bigr)$.
The remaining  coordinate changes are
$(x{:}y)\mapsto  (ax+by : cx+dy)$  where  $ad-bc\neq 0$. 


We have proved the following.

\begin{lem}  The set  of all smooth, projective 
 curves of genus 2 can be viewed as
\begin{enumerate}
\item   $\bigl((\p^1)^6\setminus\mbox{diagonals}\bigr)/\bigl(S_6\times \PGL_2\bigr)$,
equivalently as
\item   $\bigl(\sym^6(\p^1)\setminus\mbox{diagonal}\bigr)/\PGL_2$. \qed
\end{enumerate}
\end{lem}

One can see that  $\bigl(\sym^6(\p^1)\setminus\mbox{diagonal}\bigr)/\PGL_2$
has a natural structure of a normal variety   and 
from now on we set $M_2:= \bigl(\sym^6(\p^1)\setminus\mbox{diagonal}\bigr)/\PGL_2$.

In compactifying $M_2$, we need to describe what happens when the points come together. Here is a typical example. For distinct $a_i$,
consider the family  
$$
f_t:=(x-ta_1y)(x-ta_2y)(x-ta_3y)(x-ta_4y)(x-a_5y)(x-a_6y).
$$
We have 6 distinct roots for $t\neq 0$ but a 4-fold root for $t=0$.
Consider next the
coordinate change $x=tx', y=y'$ and divide by $t^4$ to obtain
$$
f'_t:=(x'-a_1y')(x'-a_2y')(x'-a_3y')(x'-a_4y')(tx'-a_5y')(tx'-a_6y').
$$
Note that $C_t\cong C'_t$ for $t\neq 0$
but $f'_0$   has only a  multiplicity 2 root, at $y'=0$. 
Variants of this trick prove the following.

\begin{lem} Let $C_t:=\bigl(z^2=f_6(x{:}y,t)\bigr)$ be a
family of genus 2 curves  whose coefficients are holomorphic functions of $t$. 
Assume that $C_t$ is smooth for $t\neq 0$. 
Then, after a coordinate change  (that depends holomorphically  on $t$)
and multiplying by a power of $t$ we get a family 
$C'_t:=\bigl(z^2=f'_6(x{:}y,t)\bigr)$ such that  
\begin{enumerate}
\item  $C_t\cong C'_t$ for $t\neq 0$ and
\item $f'_6(x{:}y,0)$ has
 either at worst double roots or 2 triple roots. \qed
\end{enumerate}
\end{lem}




\begin{cor} \label{M2.GIT.cor}
$M_2$ has a compactification $\bar M_2^{\rm GIT}$
(this is  Mumford's  GIT compactification)  whose points 
 correspond to the curves of the following types:
\begin{enumerate}
\item smooth genus 2 curves,
\item elliptic curves with 1 node,
\item rational curves with 2 nodes,
\item rational curve with a non-planar triple point,
\item rational curve with 2 cusps.  
\end{enumerate}
\end{cor}

Proof. If $f_6$ has 0,1 or 2 double roots then 
$\bigl(z^2=f_6(x{:}y)\bigr)$ is  smooth of genus 2,
 elliptic  with 1 node,
or rational with 2 nodes. 

If there are 3 double roots then we get 
  $C:=\bigl(z^2=x^2(x-y)^2(x+y)^2\bigr)$. 
This is reducible with irreducible components
$C_{\pm}:=\bigl(z=\pm x(x-y)(x+y)\bigr)$.  Choosing either of $C_{\pm}$ and 
identifying the 3 intersection points with each other gives 
a rational curve with a non-planar triple point. (There is a unique such curve.)

The last case is  $\bigl(z^2=x^3y^3\bigr)$. \qed
 \medskip
 
This is the end of the old style story: we have a compactification whose points
naturally correspond to smooth genus 2 curves and some of their degenerations.
(See \cite{MR0114819, MR2166084} for detailed descriptions of $\bar M_2^{\rm GIT}$ and  references.)
However,  from the modern point of view,  $\bar M_2^{\rm GIT}$ is  a very  unpleasant  compactification since it does not give a description 
of {\em flat families} of singular genus 2 curves.

The really big problem appears when we try to deform
the curve  $z^2=x^3y^3$. The universal deformation space is
$$
z^2=\bigl(x^3+uxy^2+vy^3\bigr)\bigl(y^3+syx^2+tx^3\bigr)
\subset \p^2(1,1,3)\times \a^4_{uvst}.
$$
If     $u=v=0$ but $st\neq 0$ (or, symmetrically, if $s=t=0$ but $uv\neq 0$) then the equations
define  disallowed curves with 1 triple root and 3 simple roots. 

We also have difficulties at the  points of type (\ref{M2.GIT.cor}.4), but these are more subtle. \qed
\medskip

The Deligne--Mumford compactification was invented to avoid these problems.

\begin{say}[Deligne--Mumford compactification  $\bar M_2$]
       The points of  $\bar M_2$
 correspond to the following curves:
\begin{enumerate}
\item smooth genus 2 curves,
\item elliptic curves with 1 node,
\item rational curves with 2 nodes,
\item 2 rational curves meeting at 3 nodes,
\item 2 elliptic curves meeting at 1 node.
\end{enumerate}

Here the cases (1--3) are the same as in (\ref{M2.GIT.cor}.1--3).
In  (\ref{M2.GIT.cor}.4)  we had $C:=\bigl(z^2=x^2(x-y)^2(x+y)^2\bigr)$.
Instead of contracting one if its irreducible components, we just keep both.

The big change is with the case (\ref{M2.GIT.cor}.5).  
A seemingly special but in fact typical   
example is given by the deformation
$$
\bigl(z^2=(x^3+t^3y^3)(y^3+t^3x^3)\bigr).
$$
For $t=0$ we get $(x^2=x^3y^3)$ and for $t^6\neq 1$ we get a smooth genus 2 curve.
Next we change the equation to
$$
z^2=t(x^3+t^3y^3)(y^3+t^3x^3).
$$
For $t\neq 0$ we did not change the fiber, but the whole family did change.
(Such changes are  actually not typical, they are only possible due to the hyperelliptic involution.)

Next we focus on the singularities at  the points  $(0{:}1,0)$ and  $(1{:}0,0)$. They are symmetric, so consider $(0{:}1,0)$. 
In affine coordinates we have 
$$
z^2=t(x^3+t^3)(1+t^3x^3).
$$
After blowing up the origin, in the new coordinates   $(x_1,  z_1, t)=(x/t, z/t, t)$ we get
$$
z_1^2=t^2(x_1^3+1)(1+t^6x_1^3).
$$
This is singular along the   line
$(z_1=t=0)$. Blow up the line  to get
$$
z_2^2=(x_1^3+1)(1+t^6x_1^3)\qtq{where} (x_1,  z_2, t)= (x_1, z_1/t, t)
$$
For $t=0$ we see the elliptic curve $(z_2^2=x_1^3+1)$. 
The other elliptic curve appears over $(1{:}0,0)$. 

Finally we can contract the birational transform of the line $(z=t=0)$ to get
a new family whose fiber over $t=0$ consists of 2 elliptic curves meeting at a node. 
\end{say}

The Deligne--Mumford compactification  $\bar M_2$ 
is a special case of the following general result, proved in
\cite{del-mum}.

 \begin{thm}[Deligne--Mumford compactification  $\bar M_g$]\label{dm.comp.thm}
$M_g$ has  a natural  compactification
whose points correspond to  {\em DM-stable curves.} 
      These are projective, connected, reduced curves $C$ 
defined by 2 properties.
 \begin{enumerate}
\item {\bf (Local)} The singularities of $C$ are   at worst nodes.
 \item  {\bf (Global)}  $\omega_C$ is ample.
  \end{enumerate}
\end{thm}







\section{Moduli spaces using Mumford's GIT}
\label{sec.4}

Mumford's GIT \cite{git} gives a recipe to construct moduli spaces of
 varieties endowed with an ample line bundle.

\begin{say}[Construction of GIT moduli spaces]\label{git.say}
The GIT approach to moduli spaces can be naturally divided into 5 steps.

\medskip
{\it Step \ref{git.say}.1.} Start with our class of  pairs   $\{(X, L)\}$.
At the beginning $X$ is any projective scheme and $L$ an ample line bundle on $X$.  Since the Hilbert function is locally constant in flat families, we may as well fix a
polynomial $p$ and let  
${\mathcal M}_p$ denote the set of all pairs $(X, L)$ where $X$ is a projective scheme, $L$ an ample line bundle on $X$  and $\chi(X, L^r)=p(r)$ for every $r$.

\medskip
{\it Step \ref{git.say}.2.} We would like to  choose an $m$ (depending on the polynomial $p$) such that $L^m$ is very ample for every 
 $(X, L)\in {\mathcal M}_p$. It turns out that there is no such $m$ that works for all $(X,L)$,  but we can choose  $m$ if we put further restrictions on $X$   (for example nodal curves, or normal surfaces, or smooth varieties, ...).
This will be a (sometimes quite hard) technical issue that does not effect the general picture. Thus for now let us denote by 
${\mathcal M}_{p,m}\subset {\mathcal M}_p$  the set of all pairs $(X, L)$ 
for which  $L^m$ is very ample and 
 $h^0(X, L^m)=\chi(X, L^m)$. 

\medskip
{\it Step \ref{git.say}.3.} Consider the images of all maps
$|L^m|: X\into \p^N$ where $N:=h^0(X, L^m)-1=p(m)-1$ and prove that they
correspond to points of a subscheme  
$$
\operatorname{EMB}_{p,m}\subset \chow(\p^N)
\qtq{or}
\operatorname{EMB}_{p,m}\subset \hilb(\p^N).
$$
(In many cases of interest these two are isomorphic, that is why I use the same notation.)
Note that $\PGL_{N+1}$ acts on $\operatorname{EMB}_{p,m} $
and the isomorphism classes of pairs $(X, L)$ correspond to the orbits of this action. (This is correct over $\c$, one needs to be more careful over non-closed fields.)

\medskip
{\it Step \ref{git.say}.4.} Mumford's GIT identifies GIT-stable and GIT-semi-stable points of the $\PGL_{N+1}$-action and constructs the corresponding coarse moduli space, denoted by   
$$
M_{p,m}^{\rm GIT}:=\operatorname{EMB}_{p,m}/\!/ \PGL_{N+1}.
$$

\medskip
{\it Step \ref{git.say}.5.}  While it is easy to define what GIT-stable and GIT-semi-stable points are, in practice  it is frequently hard to decide which points are such. Understanding GIT-stability in terms of the geometry of $(X, L)$ has been the hardest part of the project.
This was done very successfully for the moduli of sheaves; see \cite[Chap.4]{huy-leh} for a survey.

 For moduli of varieties Mumford and later Gieseker  had several major  positive results.
\begin{itemize}
\item $M_g$, the moduli of smooth curves,    Mumford  \cite{git},
\item $A_g$, the moduli of abelian varieties,    Mumford  \cite{git},
 \item $\bar M_g$, the moduli of stable curves,    Gieseker--Mumford  \cite{mumf-MR0450272, MR691308},
\item the moduli of canonical models of surfaces,  Gieseker  \cite{MR0498596}.
\end{itemize}
\end{say}

In these  results the choice of $m$ in (\ref{git.say}.2) plays a minor role. For example, for smooth curves $C$ we always choose $L:=\omega_C$.   Thus $p(t)=(2g-2)t+1-g$ and  $L^m$ is very ample for $m\geq 3$.
The hardest part is to prove that 
 the corresponding point  $[C]\in \operatorname{EMB}_{p,m}$ is GIT-stable
for every $m\geq 3$. In all the other cases, every $m\gg 1$ gives the same answer. This led Mumford to the following definition \cite[1.17]{mumf-MR0450272}.

\begin{defn} A  pair  $(X, L)$ is called
{\it asymptotically stable}  (resp.\ {\it asymptotically semi-stable})
if the corresponding point  $[X,L]\in \operatorname{EMB}_{p,m} $
is  GIT-stable  (resp.\ GIT-semi-stable) for all  $m\gg 1$. 

The following weaker variant is probably not   a useful notion in general, but it makes
the next result stronger.
$(X, L)$ is called
{\it weakly asymptotically (semi)stable}  
if   $[X,L]\in \operatorname{EMB}_{p,m} $
is  GIT-stable  (resp.\ GIT-semi-stable) for
infinitely many values of $m$. 
\end{defn}

A combination of 
\cite{MR1669716, MR1916953,  MR2103309} proves  that
for a smooth projective variety $X$ with ample canonical class $K_X$,
the pair  $(X, K_X)$ is asymptotically stable. Thus GIT  constructs the moduli space of smooth canonical models. However, we run into difficulties with
sigular ones.
This is a consequence of the next  result of   Xiaowei Wang  and 
Chenyang Xu   that  compares asymptotic stability
with KSB stability (see Definitions~\ref{KSB.stab.defn} and \ref{KSB.stab.defn.2}).

\begin{thm}\cite{xu-wang}
 Let  $0\in B$ be a smooth, pointed  curve, $B^0=B\setminus\{0\}$
and  $f^0:Y^0\to B^0$ a family of canonical models.
Assume that it has 2 compactifications
$$
\begin{array}{ccl}
Y^0 & \subset & Y^{\rm KSB}\\
f^0\downarrow\hphantom{f^0} && \ \downarrow f^{\rm KSB}\\
B^0 & \subset & B
\end{array}
\qtq{and}
\begin{array}{ccl}
Y^0 & \subset & Y^{\rm GIT}\\
f^0\downarrow\hphantom{f^0} && \ \downarrow f^{\rm GIT}\\
B^0 & \subset & B
\end{array}
$$
where $f^{\rm KSB}:Y^{\rm KSB}\to B$ is KSB stable (as in Definition~\ref{KSB.stab.defn.2}) and
$Y^{\rm GIT}_0 $ is weakly asymptotically semi-stable.

Then $Y^{\rm KSB}=Y^{\rm GIT}$.
\end{thm}

At first this would seem great news: the GIT and the KSB  compactifications seem to coincide.
However, Mumford proved in \cite[Sec.3]{mumf-MR0450272} that
 if $X$ has a point of multiplicity $>(\dim X+1)!$ then
$(X, L)$ is never asymptotically semi-stable.
There are   KSB-stable surfaces, and also canonical 3-folds, with singular points of arbitrary high multiplicity. 
Hence these are not  asymptotically semi-stable. Thus we obtain that GIT semi-stability depends too much  on the choice of $m$.

\begin{cor} \cite{xu-wang}  When we  compactify the Gieseker moduli space
of canonical models of surfaces using the $m$-canonical embeddings, 
we get infinitely many different notions of GIT-semi-stable limits, depending on $m$. \qed
\end{cor}

This already happens in rather concrete situations.

\begin{exmp}[Smooth surfaces in $\p^3$]\cite[Sec.4.2]{xu-wang}. 
Let $MS_d$ denote the moduli space of smooth surfaces of degree $d$ in  
$\p^3$. For any $m\geq 1$, using $|mH|$ we get a GIT  compactification
$$
MS_d\subset \overline{MS}_{d,m}^{\rm GIT}.
$$
The computations in \cite{xu-wang} show that, 
as $m$ varies, we get infinitely many different
notions of GIT-semi-stable limits for  $d\geq 30$. 
(Most likely the same holds for all $d\geq 5$.) 
\end{exmp}



\begin{say}[Polarization of $\operatorname{EMB}_{p,m} $]
In (\ref{git.say}.4)  we actually have one more choice to make: an ample  line bundle
on $\operatorname{EMB}_{p,m} $. Both $\chow(\p^N)$ and $\hilb(\p^N)$ naturally come with an ample
line bundle, and the traditional GIT approach uses these line bundles.

 However, as Mumford pointed out \cite[p.41]{git}, GIT does construct every projective moduli space, once we choose the correct ample  line bundle on $\operatorname{EMB}_{p,m} $. 

 The works of \cite{MR3779955, ko-pa-proj} tell us what this ample  line bundle should be for KSB-stable varieties. It would be interesting to see if a direct GIT approach works using these. 
Earlier Viehweg \cite{vieh-book} proved that the moduli space of 
 higher dimensional canonical models is quasi-projective.
\end{say}

     \section{Canonical models in higher dimensions}
\label{sec.5}

{\bf Assumption.} From now on we work over $\c$. Although everything is expected to hold over $\spec\z$, currently the
proofs are complete in characteristic 0 only.
\medskip

What are the correct higher dimensional analogs of 
 smooth, projective curves of genus $\geq 2$?
In dimension 2,  the classical answer may have been:  surfaces of general type with ample canonical class. It was also understood since the works of Du~Val \cite{DuVal34} that one could allow certain singularities. The final statement was proved by Mumford   \cite{MR0141668b}: 
If $S$ is a smooth, projective surface of general type then
its canonical ring   
$$
R(S, \omega_S):=\tsum_{m\geq 0} H^0\bigl(S, \omega_S^m\bigr)
$$
is finitely generated and $S^{\rm can}:=\proj R(S, \omega_S)$---called the  {\it canonical model} of $S$---is a
surface with ample canonical class and at worst Du~Val singularities.

The following was proved by Mori in dimension 3 \cite{mori-MR924704} 
(with important contributions by  Kawamata, Koll\'ar, Reid,  Shokurov)
and by
Hacon and McKernan  in higher dimensions
(with important contributions by  Birkar, Cascini, Corti,  Shokurov).

\begin{thm}[Canonical models] 
Let $X$ be a smooth, projective variety of general type. Then
its canonical ring   
$$
R(X, \omega_X):=\tsum_{m\geq 0} H^0\bigl(X, \omega_X^m\bigr)
$$
is finitely generated and $X^{\rm can}:=\proj R(X, \omega_X)$---called the  {\it canonical model} of $X$---is a
variety with ample canonical class and at worst ``canonical''  singularities.
\end{thm}

We have not yet defined what these   ``canonical''  singularities are.
In fact historically this was a major difficulty that was resolved by Reid  \cite{r-c3f}.  First we need to recall what
the canonical sheaf  of a singular variety is.

\begin{say}[Canonical class or  sheaf I]\label{cc.ds.defn}
First let $U$ be a smooth variety.
The determinant of its cotangent bundle
$\det \Omega_U^1$ is a line bundle, frequently denoted by $\Omega_U^n$  (where $n=\dim U$), though algebraic geometers prefer the notation $\omega_U$ and call it the
{\it canonical  bundle} of $U$.  The divisor of a rational section of  $\omega_U$ is called a {\it canonical  divisor} on $U$. The corresponding linear equivalence class is denoted by $K_U$ and called the
{\it canonical  class} of $U$. (For reasons that are not important for us now,
 $\omega_U$ is frequently called  the
{\it dualizing  bundle or sheaf} of $U$.) 

Let $X$ be a normal variety and       $X^0\subset X$ its 
 smooth locus.  Note that  $X\setminus X^0$ has codimension $\geq 2$ in $X$.
The closure of a canonical divisor  $K_{X^0}$ defines a 
canonical divisor  $K_{X}$ and the
push-forward of the canonical bundle  $\omega_{X^0}$ defines the 
canonical sheaf $\omega_{X}$. 

In general  $\omega_{X}$ is not locally free, so its tensor powers can be very complicated. It is better to use instead the reflexive powers
$\omega_X^{[m]}$, defined in any of the following ways:
\begin{enumerate}
\item the push-forward of  $\omega_{X^0}^m$,
\item the double dual of  $\omega_X^{\otimes m}$ or
\item the sheaf  $\o_X(mK_X)$.
\end{enumerate}
\end{say}

\begin{exmp} \label{omega.exmp}
The following examples are fairly typical.

(\ref{omega.exmp}.1) Let  $X:=(g=0)\subset \c^n$ be a hypersurface.
As discussed in \cite[Sec.III.6.4]{shaf}, $\omega_X$ is locally free and a generator of $\omega_X$ can be written as
$$
\sigma_i:=(-1)^i\frac{dx_1\wedge\dots \wedge \widehat{dx_i}\wedge\dots \wedge dx_n}{\partial g/\partial x_i}.
$$
Despite appearances, the $\sigma_i$ glue to a regular section $\sigma_X$ of
$\omega_X$. However,  
a pull-back of $\sigma$ frequently has poles since we are dividing by ${\partial g/\partial x_i} $.

(\ref{omega.exmp}.2) Let  $Y:=\c^n/G$ be a quotient singularity,  where $G$ is a finite subgroup of $\GL_n$. Assume for simplicity that $G$-action is free outside a codimension $\geq 2$ subset of $\c^n$.
In general $\omega_Y$ is not locally free but, for $m:=|G/G\cap \SL_n|$,
$$
(dx_1\wedge\dots \wedge dx_n)^{\otimes m}
$$
descends to   a generator of $\omega_Y^{[m]}$, which is thus locally free.
\end{exmp}

 \begin{defn}[Canonical singularities] Informally, we say that a normal variety $X$ has canonical singularities if one can pull back its canonical class/sheaf the same way as for smooth varieties.

Formally, for  every  resolution   $p:Y\to X$   we assume the following.
\begin{enumerate}
\item (Pull-back of forms.) For every $m\geq 0$ the  pull-back of forms on the smooth locus extends to  $p^*\omega_X^{[m]}\to \omega_Y^{m}$.
(It is enough the check this condition  for one resolution.)

\item (Pull-back of canonical divisors.)  $p^*K_X$ makes sense. Note that pull-back is not defined for divisors but it is defined for Cartier divisors. 
More generally, if $mK_X$ is Cartier for some $m>0$ then
$p^*K_X$ makes sense as a rational linear combination of prime divisors on $Y$.
\end{enumerate}

By accident, if (1) holds for a surface $S$  then $K_S$ is Cartier.
Thus it took some time to realize that in higher dimensions we need to add (2) as a new assumption and that only a multiple of $K_X$ needs to be Cartier.
  \end{defn}

We can now give an 
 internal definition of canonical models.
We declare these to be the   correct higher dimensional analogs of 
 smooth, projective curves of genus $\geq 2$.

\begin{defn}[Canonical model] A 
normal, projective variety $X$ is a canonical model iff 
\begin{enumerate}
\item {\bf (Local)}   $X$ has canonical singularities and
\item {\bf (Global)}   $\omega_X$ is ample, that is, $\omega_X^{[m]}$ is locally free and ample for some $m>0$.
 \end{enumerate}
\end{defn}

\section{What is a node?}
\label{sec.6}

DM-stable curves have nodes, thus it is reasonable to expect that
in general we need to allow
 higher dimensional generalizations of nodes in order to get a compact moduli space. The following 3 characterizations of nodes were well known.
The problem was that nodes have many other special properties and it was not clear which ones should  be generalized.

\begin{prop}[Characterizations of nodes]\label{char.nodes}
Let $(0\in C)$ be a reduced curve singularity such that $\omega_C$ is free with
generating section $\sigma$. Then $(0\in C)$ is a node iff any of the following holds.
\begin{enumerate}
\item ({\it Using resolution}) 
Let $p:\bar C\to C$ denote the normalization. 
Then   $p^*\sigma$ has only
 simple poles, as a rational section of $\omega_{\bar C}$.
\item ({\it Using  volume})  Let $0\in U\subset C$ be a relatively compact open subset. Although the  volume 
$\sqrt{-1} \int_{U}  \sigma\wedge \bar \sigma$ is  $\infty$,
it has only logarithmic growth, that is, 
$$
\sqrt{-1}\int_U  |g|^{\epsilon} \cdot \sigma\wedge \bar \sigma <\infty
$$
for every $g$ vanishing at $0\in C$ and $\epsilon>0$.
\end{enumerate}
\noindent If $C$ is smoothable then these  are also equivalent to the following.
\begin{enumerate}\setcounter{enumi}{2}
\item ({\it Using deformations}) For every 1-parameter smoothing $p:{\mathbf C}\to \Delta$ of $C\cong {\mathbf C}_0$ over the disc $\Delta$, the total space ${\mathbf C}$ has
canonical singularities.
\end{enumerate}
\end{prop}

\begin{say}[Comments] We can write a node as
$C:=(xy=0)\subset \c^2$. By (\ref{omega.exmp}.1),  $\omega_C$ has a 
generating section $\sigma$ such that 
$$
\sigma=\tfrac{dx}{x} \mbox{ on the $x$-axis, }\quad  \sigma=-\tfrac{dy}{y} \mbox{ on the $y$-axis.}
$$
Thus $\sigma$ has only simple poles one the resolution, which is the disjoint union of the 2 axes. 
 (The sign is not important for us, in general the sum of the 2 residues must be zero.)
By contrast let $D:=(x^2=y^3)\subset \c^2$ be a cusp.
By (\ref{omega.exmp}.1) a generating section of $\omega_D$ is  $\sigma_D:=\tfrac{dx}{3y^2}=-\tfrac{dy}{2x}$. The normalization is
$p:\c\to D$ given by $t\mapsto (t^3, t^2)$. Thus
$p^*\sigma_D=\tfrac{dt^3}{3t^4}=\tfrac{dt}{t^2}$ has double pole.

Condition (\ref{char.nodes}.2) 
 boils down to the elementary local formulas
$$
\sqrt{-1}
\int_{|x|\leq 1}  \frac{dx}{x}\wedge \frac{d\bar x}{\bar x} =\infty,
\qtq{but}
\sqrt{-1}
\int_{|x|\leq 1}  |x|^{\epsilon}\frac{dx}{x}\wedge \frac{d\bar x}{\bar x} <\infty.
$$
As for (\ref{char.nodes}.3), every smoothing of a node 
is of the form  $(xy=t^n)$ for some $n$ and these are all  canonical.
 For 
the cusp  $(x^2+y^3=0)$ the  deformations  $(x^2+y^3=t^n)$ are   canonical only for $n\leq 5$.

The caveat before (\ref{char.nodes}.3) is necessary since 
Mumford showed that not all curve singularities are smoothable \cite{MR0460338}.
  \end{say}

\section{Compactification using minimal models}
     \label{sec.7} 
   
The original construction of \cite{del-mum} used the Jacobian to construct stable limits. Another proof, using  Mumford's semi-stable reduction theorem \cite{MR0335518} first constructs a reduced degeneration with only nodes and then uses the theory of minimal models of surfaces to contract the  unnecessary rational components. The key idea of \cite{ksb} is to use this second approach of Mumford
and  the following  theorem on  relative canonical models, which
 is a special case of the general results of Mori, Hacon and McKernan, to be discussed in Paragraph~\ref{hist.of.pf}.

\begin{thm} Let  $0\in B$ be a smooth, pointed  curve, $B^0=B\setminus\{0\}$
and  $f^0:Y^0\to B^0$ a family of canonical models.
Then there is a unique extension to
$$
\begin{array}{ccc}
Y^0 & \subset & Y\\
f^0\downarrow\hphantom{f^0} && \hphantom{f}\downarrow f\\
B^0 & \subset & B
\end{array}
$$
such that  $Y$ has canonical singularities and
 $\omega_Y$  is ample on every fiber of $f$.
\end{thm}

\cite{ksb}  used this to give an abstract characterization of the resulting central fiber $Y_0$ in the semi-stable reduction case; this is recalled in
Theorem~\ref{slc.defn.2}. 

For now  we turn history around and first define the notion of  semi-log-canonical singularities as a direct generalization of the characterizations of nodes given in Proposition~\ref{char.nodes}.

 Since DM-stable curves are not normal,  we also need a non-normal version of Definition~\ref{cc.ds.defn}.

\begin{say}[Canonical class or  sheaf II]\label{cc.ds.defn.2}
 Recall that, by Serre's theorem, a variety is normal iff it is smooth at codimension 1 points and satisfies the condition $S_2$ (= Hartogs's extension theorem). 
As a slight generalization, we say that $X$ is  {\it demi-normal} iff it is either smooth or has nodes at codimension 1 points and satisfies the condition $S_2$. 

If this holds then there is an open subset  $X^0\subset X$  such that
$\omega_{X^0}$ is locally free and 
$X\setminus X^0$ has codimension $\geq 2$ in $X$.
We can thus define  the reflexive powers
$\omega_X^{[m]}$ in either of the following ways:
\begin{enumerate}
\item the push-forward of  $\omega_{X^0}^m$,
\item the double dual of  $\omega_X^{\otimes m}$.
\end{enumerate}


\end{say}

 \begin{defn}[Semi-log-canonical singularities]\label{slc.defn}
 Let  $X$ be a demi-normal variety such that $\omega_X^{[m]}$ is locally free  for some $m>0$.  
We say that $X$ has {\it semi-log-canonical singularities} iff the following equivalent conditions hold.
\begin{enumerate}
\item ({\it Resolutions I}) Let    $f:Y\to X$ be a resolution and $E\subset Y$ the  reduced exceptional divisor $E$. Then the pull-back map on forms extends to 
$$
f^*\bigl(\omega_X^{[m]}\bigr)\to \omega_Y^{m}(mE).
$$
\item ({\it Resolutions II}) Let    $f:Y\to X$ be a resolution and $E\subset Y$ the  reduced exceptional divisor $E$. Then the pull-back map on forms extends to 
$$
f^*\bigl(\omega_X^{[r]}\bigr)\to \omega_Y^{r}(rE) \quad \forall r.
$$
\item ({\it Local volume})  Let $U\subset X$ be any relatively compact open
and     $\sigma^m$ a  generating section of 
$\omega_X^{[m]}|_{\bar U}$. Then 
$$
i^{n(2-n)}\int_U  |g|^{\epsilon} \cdot \sigma\wedge \bar \sigma <\infty
$$
for every polynomial $g$ vanishing on $\sing U$ and $\epsilon>0$.
\end{enumerate}
The above conditions are analogs of (\ref{char.nodes}.1--2). There is also a generalization of 
(\ref{char.nodes}.3) using deformations, but it needs an extra condition
that we explain later; see Theorem~\ref{slc.defn.2}.

If $X$ has a  1-parameter smoothing $p:{\mathbf X}\to \Delta$ such that  $\omega_{\mathbf X}^{[M]}$ is locally free for some $M>0$,  then (1--3)  are also equivalent to the following. 
\begin{enumerate}\setcounter{enumi}{3}
\item ({\it Using deformations})  For every such smoothing, the total space ${\mathbf X}$ has
canonical singularities.
\end{enumerate}
\end{defn}

We can now define the higher dimensional generalizations of DM-stable curves.

\begin{defn}[KSB-stable variety]\label{KSB.stab.defn}
A projective variety $X$ is {\it KSB-stable} if the following hold.
\begin{enumerate}
\item ({\bf Local})  $X$ has only semi-log-canonical singularities and
\item ({\bf Global})  $\omega_X$ is ample.
\end{enumerate}
\end{defn}

\begin{say}[Comments on \ref{slc.defn}]
The equivalence of conditions (\ref{slc.defn}.1--3) is quite straightforward
but the equivalence of these with  (\ref{slc.defn}.4) is  hard; see
Theorem~\ref{slc.defn.2} for an explanation and an outline of the proof.

In (\ref{slc.defn}.3) note that although $\sigma$ is a multiple-valued section of $\omega_X|_U$,
the ambiguity is an $m$th root of unity, so 
$\sigma\wedge \bar \sigma $ is well defined. The power of $\sqrt{-1}$ needed depends on
the orientation conventions.
\end{say}

\section{Stable families and their moduli space}
     \label{sec.8}

Once we settled on stable curves as the objects of our moduli theory,
we conclude that a  ``good family'' of curves   is a
flat, projective  morphisms with stable fibers.
However, in higher dimensions one more twist awaits us. 

\begin{surprise}  Flat, projective  morphisms with stable fibers
do {\bf not} give the ``good families'' of a moduli theory.
\end{surprise}

We illustrate this with an example.

\begin{exmp}
      Consider the flat family of varieties in $\p^5_{\mathbf x}\times \a^2_{st}$ given by the equations 
$$
X:=\left(\rank
\left(
\begin{array}{ccc}
x_0 & x_1 & x_2\\
x_1+sx_4 & x_2+tx_5 & x_3
\end{array}
\right)
\leq 1\right)
\subset \p^5_{\mathbf x}\times \a^2_{st}.
$$
 We claim that  $X_{st}$ has semi-log-canonical  singularities iff 
either $(s,t)=(0,0)$ or  $st\neq 0$.  Thus having
semi-log-canonical  singularities is not a locally closed
condition. We modify this example below to show  that 
being KSB-stable is also  not a locally closed condition in flat families of normal, projective varieties.

Depending on whether $s,t$ are zero or not, linear coordinate changes 
show that we have only 3 types of fibers as follows.

Case 1: $st\neq 0$. Then $X_{st}\cong \p^1\times\p^2$ using 
$$
\left(
\begin{array}{ccc}
x_0 & x_1 & x_2\\
x_1+sx_4 & x_2+tx_5 & x_3
\end{array}
\right)
\to 
\left(
\begin{array}{ccc}
x_0 & x_1 & x_2\\
x_4 & x_5 & x_3
\end{array}
\right)
$$
 
Case 2: $s=t= 0$. Then $X_{st}\cong$  cone over $\p^1\hookrightarrow\p^3$ using
$$
\left(
\begin{array}{ccc}
x_0 & x_1 & x_2\\
x_1+sx_4 & x_2+tx_5 & x_3
\end{array}
\right)
\to 
\left(
\begin{array}{ccc}
x_0 & x_1 & x_2\\
x_1 & x_2 & x_3
\end{array}
\right)
$$
Thus the singularities of $X_{st}$ are of the following form:
the quotient singularity
$\c^2/\z_3$  (using the diagonal action) times a smooth factor. 
 
Case 3: $s=0, t\neq  0$.  Then $X_{st}\cong$  cone over $B_0\p^2\cong \f_1\hookrightarrow\p^4$ using.
$$
\left(
\begin{array}{ccc}
x_0 & x_1 & x_2\\
x_1+sx_4 & x_2+tx_5 & x_3
\end{array}
\right)
\to 
\left(
\begin{array}{ccc}
x_0 & x_1 & x_2\\
x_1 & x_5 & x_3
\end{array}
\right)
$$
Since the hyperplane class of $\f_1\hookrightarrow\p^4$ is not proportional to the canonical class of $\f_1$, thus no multiple of the   canonical class of $X_{st}$ is Cartier. 

We can also compute that
$$
\bigl(K_{X_{00}}^3\bigr)=-56\tfrac{8}{9}  \qtq{and} \bigl(K_{X_{st}}^3\bigr)=-54
\qtq{for} st\neq 0,
$$
thus the volume of the fibers jumps in the 1-parameter family
$$
Y:=\left(\rank
\left(
\begin{array}{ccc}
x_0 & x_1 & x_2\\
x_1+tx_4 & x_2+tx_5 & x_3
\end{array}
\right)
\leq 1\right)
\subset \p^5_{\mathbf x}\times \a^1_{t}.
$$
In Cases 1 and 2 the $X_{st}$ are Fano varieties, hence their canonical class is not at all ample. We can however easily get KSB-stable examples from these by
taking a cyclic cover. To be very concrete,
for $m\geq 1$ let $Y_m\subset \p^6_{\mathbf x}\times \a^1_{t}$ be the family  defined by the equations
$$
\tsum x_i^m=0\qtq{and}
\rank
\left(
\begin{array}{ccc}
x_0 & x_1 & x_2\\
x_1+tx_4 & x_2+tx_5 & x_3
\end{array}
\right)
\leq 1.
$$
An easy computation shows that 
$Y_m\to \a^1_t$ is  a flat, projective morphism with stable  fibers
for $m\geq 5$,
 yet the volume of the fibers is not locally constant.
 \end{exmp}

The following result connects global and local properties of deformations.

\begin{thm}\cite{k-howmuch}  \label{vol.stab.thm}
Let $X\to S$ be a  flat, projective  morphisms with KSB-stable fibers of dimension $n$ and  $S$ reduced. The following are equivalent.
\begin{enumerate}
\item  The volume of the fibers $s\mapsto \bigl(K_{X_s}^n\bigr)$ is locally constant.
\item The plurigenera  $s\mapsto h^0\bigl(X_s, \omega_{X_s}^{[m]}\bigr)$ are  locally constant for every $m$.
\item $\omega_{X/S}^{[m]}$ is flat and commutes with base change  $\forall m$.
\end{enumerate}
\end{thm}

The   volume  $\bigl(K_{X}^n\bigr)$ is the most basic numerical invariant of a stable variety, and we would definitely like it to be locally constant in 
``good'' families. The deformation invariance of all the  
 plurigenera  $P_m(X):=h^0\bigl(X, \omega_{X}^{[m]}\bigr)$ is an added bonus.
Condition (\ref{vol.stab.thm}.3) seems complicated, but 
\cite{k-hh} proves that it defines a good moduli functor over arbitrary schemes. we can thus choose it to define the notion of sable families in general.

  \begin{defn}[KSB-stable families]\label{KSB.stab.defn.2}
      Let  $S$  be an arbitrary  scheme. A  morphism $f:X\to S$ is 
{\it KSB-stable}     iff
\begin{enumerate}
\item  $f$ is flat, projective  with KSB-stable fibers and
\item $\omega_{X/S}^{[m]}$ is flat and commutes with base change  $\forall m$.
\end{enumerate}
\end{defn}

We can now state the main theorem for the higher dimensional analogs
of $\bar M_g$.

\begin{thm}[Moduli of stable varieties]
Fix  positive $n\in \n, d\in \q$. 
Then the functor  $\bar{\mathcal M}_{n.d}$ of KSB-stable families of dimension $n$ and volume $d$ has a projective  coarse moduli space $\bar M_{n,d}$. 
\end{thm}

The  moduli properties of $\bar M_{n,d}$ are as good as for $\bar M_g$, 
but as a scheme  $\bar M_{n,d}$ can be  much more complicated and very few examples have been computed explicitly \cite{ale-abvar, MR2956036}.

   \begin{say}[History of the proof]\label{hist.of.pf}
    
For surfaces the background material needed Mori's minimal model program  for 3-folds; see
\cite{km-book} for an introduction. 

Koll\'ar--Shepherd-Barron proposed to use stable surfaces for compactifications \cite{ksb}. The paper gave a classifications of semi-log-canonical surface singularities and proved that  $\bar M_{2,d}$ exists, is locally of finite type and satisfies the valuative criterion of properness.
Alexeev \cite{al94} then proved that the $\bar M_{2,d}$ is
of finite type, hence proper.  
Earlier  Koll\'ar proved that properness implies projectivity \cite{k-proj}.

In higher dimensions we need the  minimal model program. While this is not yet fully known, the works of 
Hacon and Xu  \cite{hacon-xu-brep, MR3032329} prove all instances of the  minimal model program that are need in moduli theory.

The basic existence result for $\bar M_{n,d}$ now follows the same path as for surfaces, except that
it turned out to be quite troublesome to deal with degenerations of reducible stable varieties. This was accomplished by the gluing theory of 
\cite{k-source}, see \cite[Chap.5]{kk-singbook} for a detailed treatment.
As before,  these results establish that $\bar M_{n,d}$ exists, is locally of finite type and satisfies the valuative criterion of properness.

 Karu \cite{karu} proved that all irreducible components of $\bar M_{n,d}$ are of finite type. The proof that $\bar M_{n,d}$ itself is of  finite type is a  hard result of
Hacon-McKernan-Xu \cite{hmx-bounded}. These then imply that $\bar M_{n,d}$ is proper.
Projectivity required substantial new ideas, it was settled by  Fujino and  Kov\'acs--Patakfalvi, see 
\cite{MR3779955, ko-pa-proj}. 
Earlier Vakil showed that the local structure of $\bar M_{n,d}$ can be  arbitrarily bad
\cite{MR2227692}.

A comprehensive treatment is planned in \cite{k-modbook}.
\end{say}

We can now state the correct analog of (\ref{char.nodes}.3),
which explains item (4) in Definition~\ref{slc.defn}  of semi-log-canonical singularities.

 \begin{thm}[Deformation of semi-log-canonical singularities]\label{slc.defn.2}
 Let  $(0\in D)$ be a demi-normal singularity such that $\omega_D^{[m]}$ is locally free  for some $m>0$.

Assume that $D$ is a Cartier divisor  $D\cong (g=0)$ on  a normal singularity $(0\in Y)$ such that 
 $\omega_X^{[M]}$ is locally free  for some $M>0$ and 
$Y\setminus D$ has canonical singularities only. 
Then $D$ has semi-log-canonical singularities iff the cyclic covers
$$
Y_r:=(t^r=g)\subset Y\times \a^1_t
$$
have canonical singularities only for every $r$. 
\end{thm}

Proof.  We identify $D$ with $(t=0)\subset Y_r$. 
Take a resolution  $f_r:Z_r\to Y_r$. Write
 $$
K_{Z_r}=f_r^*K_{Y_r}+J \qtq{and} f_r^*D=D_{Y_r}+E_r.
$$
By  Mumford's semi-stable reduction theorem \cite{MR0335518}
we  may choose $r$ such  that all  coefficients of $E_r$ equal 1.
Now replace $Y$ by $Y_r$ and drop the subscript $r$ from the notation.
By the adjunction formula, 
$$
K_{D_Y}= \bigl(K_Y+D_Y\bigr)|_{D_Y}=\bigl(f^*(K_X+D)+J-E\bigr)|_{D_Y}
= f^*K_D+(J-E)|_{D_Y}.
$$
If $Y$ has canonical singularities then $J$ is effective, so every divisor appears in 
$(J-E)|_{D_Y} $ with coefficient $\geq -1$. Thus $D$ has
semi-log-canonical singularities using the characterization
(\ref{slc.defn}.1). 

Conversely, if $D$ has
semi-log-canonical singularities then  every divisor appears in 
$(J-E)|_{D_Y} $ with coefficient $\geq -1$, hence 
 every divisor appears in 
$J|_{D_Y} $ with coefficient $\geq 0$. Unfortunately this does not seem to give
any information on those irreducible components of $J$ that are disjoint from $D_Y$. 
 
However  Shokurov conjectured that this is not a problem.
There are 2 ways to approach this. By a judicious choice of $Y$ we might arrange that the irreducible component in question meets $D_Y$, or one can prove that the coefficients of $J$ have certain convexity properties.
Both of these are quite interesting, see
\cite[Sec.17]{k-etal},   Kawakita \cite{kawakita},   de~Fernex-Koll\'ar-Xu \cite{dkx}  or  \cite[Sec.4.1]{kk-singbook}.

\def\cprime{$'$} \def\cprime{$'$} \def\cprime{$'$} \def\cprime{$'$}
  \def\cprime{$'$} \def\dbar{\leavevmode\hbox to 0pt{\hskip.2ex
  \accent"16\hss}d} \def\cprime{$'$} \def\cprime{$'$}
  \def\polhk#1{\setbox0=\hbox{#1}{\ooalign{\hidewidth
  \lower1.5ex\hbox{`}\hidewidth\crcr\unhbox0}}} \def\cprime{$'$}
  \def\cprime{$'$} \def\cprime{$'$} \def\cprime{$'$}
  \def\polhk#1{\setbox0=\hbox{#1}{\ooalign{\hidewidth
  \lower1.5ex\hbox{`}\hidewidth\crcr\unhbox0}}} \def\cdprime{$''$}
  \def\cprime{$'$} \def\cprime{$'$} \def\cprime{$'$} \def\cprime{$'$}
\providecommand{\bysame}{\leavevmode\hbox to3em{\hrulefill}\thinspace}
\providecommand{\MR}{\relax\ifhmode\unskip\space\fi MR }
\providecommand{\MRhref}[2]{%
  \href{http://www.ams.org/mathscinet-getitem?mr=#1}{#2}
}
\providecommand{\href}[2]{#2}

\bigskip

\noindent  Princeton University, Princeton NJ 08544-1000

{\begin{verbatim} kollar@math.princeton.edu\end{verbatim}}

\end{document}